\documentclass[a4paper,10pt]{amsart}
\usepackage{pb-diagram}
\usepackage{xypic}
\usepackage{amscd}

\usepackage{bm}
\usepackage{amsmath,amssymb,amscd,bbm,amsthm,mathrsfs,hyperref}
\usepackage{latexsym}
\usepackage{amsfonts}
\usepackage{amsthm}
\usepackage{indentfirst}
 \newtheorem{thm}{Theorem}[section]
 
 \newtheorem{lem}[thm]{Lemma}
 \newtheorem{prop}[thm]{Proposition}
 \newtheorem{defn}[thm]{Definition}
 \newtheorem{ex}[thm]{Example}
 \newtheorem{rem}[thm]{Remark}

\newtheorem{quest}[thm]{Question}
\def\k{\mathbbm{k}}

 \newcommand{\Hom}{\mathrm{Hom}}

\title{A note on homologically smooth connected cochain DG algebras}
\author{X.-F. Mao}
\address{Department of Mathematics, Shanghai University, Shanghai 200444, China}
\email{xuefengmao@shu.edu.cn}

\author{J.-F.Xie}
\address{Department of Mathematics, Shanghai University, Shanghai 200444, China}
\email{jianfengxie@yahoo.com}
\date{}

\begin{document}
 \def\abstactname{abstract}
\begin{abstract}
In this paper, we obtain two interesting results on homologically smooth connected cochain DG algebras.
More precisely, we show that any Koszul DG module in $\mathrm{D_{fg}}(A)$ is compact, when $A$ is a homologically smooth connected cochain DG algebra with a Noetherian cohomology graded algebra $H(A)$.  And we prove that the homologically smoothness of $A$ is equivalent to  $$\mathrm{D_{fg}}(A)=\mathrm{D}^c(A),$$ if $A$ is a Koszul connected cochain DG  algebra such that $H(A)$ is a Noetherian graded algebra with a balanced dualizing complex.
\end{abstract}


%16E45 Differential graded algebras and applications
%16E10 Homological dimension

%%% ----------------------------------------------------------------------
\maketitle
%%% ----------------------------------------------------------------------

\section*{introduction}
Throughout this paper, $\k$ is a algebraically closed field with $\mathrm{char}\k=0$, and $A$ is a connected cochain DG algebra over $\k$. One sees that its cohomology graded algebra
  $$H(A)=\bigoplus_{i=0}^{+\infty}\frac{\mathrm{ker}(\partial_{A}^i)}{\mathrm{im}(\partial_{A}^{i-1})}$$
 is a connected graded algebra.
 For any cocycle element $z\in \mathrm{ker}(\partial_{\mathcal{A}}^i)$, we write $\lceil z \rceil$ as the cohomology class in $H(\mathcal{A})$ represented by $z$.  We assume that the reader is familiar with basics on DG homological algebra. If this is not the case,  we refer
to \cite{AFH,FJ2,MW1,MW2} for more details on them.

 The derived category of DG $A$-modules is denoted by $\mathrm{D}(A)$. A DG $A$-module $M$ is called \emph{cohomologically finitely generated} if $H(M)$ is a finitely generated $H(A)$-module. The full triangulated subcategory of $\mathrm{D}(A)$
consisting of cohomologically finitely generated DG $A$-modules is denoted by $\mathrm{D_{fg}}(A)$. Similarly, a DG $A$-module $M$ is called  \emph{cohomologically locally finite} if each $H^i(M)$ is finite dimensional. The full
subcategory of $\mathrm{D}(A)$ consisting of cohomologically locally finite DG $A$-modules is denoted by $\mathrm{D_{lf}}(A)$. Let $\mathrm{D^{+}}(A)$, $\mathrm{D^{-}}(A)$ and
$\mathrm{D^{b}}(A)$ be the full subcategories of $\mathrm{D}(A)$,
whose objects are \emph{cohomologically bounded below}, \emph{cohomologically bounded above} and
\emph{cohomologically bounded}, respectively.

According to the definition of compact objects (cf. \cite{Kr1,Kr2}) in a triangulated category with
arbitrary coproduct, a DG $A$-module is called \emph{compact}, if the functor $\Hom_{\mathrm{D}(A)}(M,-)$ preserves
all coproducts in $\mathrm{D}(A)$. Since $A$ is
an augmented DG algebra, a DG $A$-module $M$ is compact, if and only if it is
in the smallest thick subcategory of $\mathrm{D}(A)$ containing ${}_AA$ (cf. \cite[Theorem 5.3]{Kel}).
To use the language of topologists, a DG A-module is compact
if it can be built finitely from ${}_AA$, using suspensions and distinguished triangles. The full
subcategory of $\mathrm{D}(A)$ consisting of compact DG $A$-modules is denoted by $\mathrm{D}^c(A)$.
Compact DG modules play the same role as finitely presented
modules of finite projective dimension do in ring theory (cf. \cite{Jor2}). By \cite[Proposition 3.3]{MW1}, a DG $A$-module $M$ is compact if and only if $M$ has a minimal semi-free resolution $F_M$, which has a finite semi-basis.

One sees easily that any compact DG $A$-module is a cohomologically finitely generated $A$-module. However, the converse is generally not true. For example, if $A$ is not homologically smooth then ${}_A\k$ is not compact while $\k\in \mathrm{D_{fg}}(A)$.  Recall that a connected cochain DG algebra $A$ is \emph{homologically smooth} if ${}_{A^e}A$ is compact, which is equivalent to $\k \in \mathrm{D}^c(A)$ (cf. \cite[Corollary 2.7]{MW3}). Let $A$ be a connected cochain DG algebra such that $H(A)$ is a Noetherian graded algebra with $\mathrm{gl.dim}H(A)<\infty$. Then by the existence of Eilenberg-Moore resolution, one sees that any cohomologically finitely generated DG $A$-module is compact, which implies that $A$ is homologically smooth. And \cite[Example 3.12]{MW2} indicates that there are homologically smooth connected cochain DG algebras whose cohomology graded algebras are Noetherian graded algebras with infinite global dimension.
So the homologically smoothness of $A$ is weaker than $\mathrm{gl.dim}H(A)<\infty$ when $H(A)$ is Noetherian. It is natural to ask the following question.
\begin{quest}\label{mainquest}
{\rm Assume that $A$ is a homologically smooth connected cochain DG algebra such that $H(A)$ is Noetherian.
Is any cohomologically finitely generated DG $A$-module compact ?}
\end{quest}

In DG context, the homologically smoothness of DG algebras is analogous to the regular property of graded algebras. The research on  this fundamental property of DG algebras have attracted many people's interests.
 In \cite{HW1},  He-Wu introduced the concept of
Koszul DG algebras, and obtained a DG version of the Koszul duality for Koszul, homologically smooth and Gorenstein DG algebras.
The author and Wu \cite{MW2} proved that any homologically
smooth connected cochain DG algebra $A$ is cohomologically unbounded
unless $A$ is quasi-isomorphic to the simple algebra $\k$. And it was
proved that the $\mathrm{Ext}$-algebra
 of a homologically smooth DG algebra $A$
is Frobenius if and only if both $\mathrm{D^b_{lf}}(A)$ and
$\mathrm{D^b_{lf}}(A\!^{op})$ admit Auslander-Reiten triangles.  In \cite{Sh}, Shklyarov developed
a Riemann-Roch Theorem for homologically smooth DG algebras.
Besides these, some important classes of DG algebras are homologically smooth. For example, Calabi-Yau DG algebras introduced by Ginzburg in \cite{Gin} are homologically smooth by definition. Especially, non-trivial Noetherian DG down-up algebras and DG polynomial algebras are  Calabi-Yau DG algebras by \cite{MHLX} and \cite{MGYC}, respectively.

The importance of homologically smooth DG algebras motivates us to consider Question \ref{mainquest}.
In this paper, we partially solve it by the following interesting results (see Theorem \ref{koszulcomp}).
\\
\begin{bfseries}
Theorem \ A.
\end{bfseries}
Let $A$ be a homologically smooth connected cochain DG algebra such that $H(A)$ is Noetherian. Then
any Koszul DG module in $\mathrm{D_{fg}}(A)$ is compact.

 Here, a DG $A$-module $M$ is called \emph{Koszul} if it has a semi-free resolution $P$, which has a semi-basis concentrated in a single degree. If ${}_A\k$ is a Koszul DG module, then we say $A$ is a Koszul DG algebra.
 \begin{rem}\label{onkoszul} {\rm Note that the definition of Koszul DG modules is different from \cite[Definition 2.1]{HW2} and \cite[Definition 5.3]{Jor3} since the single degree here can be non-zero. By the definition, one sees that any Koszul DG $A$-module is an object in $\mathrm{D}^+(A)$. Suppose that $M$ is a Koszul DG $A$-module in $\mathrm{D}^+(A)$. Then by  \cite[Proposition 2.4]{MW1}, one sees that $M$ admits a minimal semi-free resolution $F_M$, which has a semi-basis concentrated in degree $\inf\{i|H^i(M)\neq 0\}$. Hence our definition of Koszul DG algebras
coincides with those of He-Wu and J$\o$rgensen in \cite{HW1} and \cite{Jor3}, respectively.}
  \end{rem}
  In \cite{Jor3}, J$\o$rgensen developed a duality between $\mathrm{D_{fg}}(A)$ and $\mathrm{D_{fg}}(A^{op})$. Under the additional condition that $H(A)$ admits a balanced dualizing complex, he obtained some interesting results on the Ext and Castelnuovo-Mumford regularities of DG modules in $\mathrm{D_{fg}}(A)$ (see \cite[Theorem 5.7]{Jor3}).  Inspired from \cite{Jor3}, we prove the following theorem (see Theorem \ref{criteria}).
\\
\begin{bfseries}
Theorem \ B.
\end{bfseries}
Let $A$ be a Koszul connected cochain DG algebra such that $H(A)$ is a
Noetherian graded algebra with a balanced dualizing complex. Then
$A$ is homologically smooth if and only if
$\mathrm{D}^c(A)=\mathrm{D_{fg}}(A)$.

In the last section, we give Example \ref{ex1} and Example \ref{ex2} to explain that there are Koszul, homologically smooth connected cochain DG algebras
whose cohomology graded algebras are Noetherian graded algebras with balanced dualizing complexes and infinite global dimension. We can apply  Theorem B to them. While the DG algebra in Example \ref{ex3} is a non-Koszul homologically smooth DG algebra. So we can only apply Theorem A to it.

%%% ----------------------------------------------------------------------
\section{Ext and Castelnuovo-Mumford regularities of DG modules}
In this section, we review the Ext and Castelnuovo-Mumford regularities of DG modules. These two invariants of DG modules were introduced and studied in \cite{Jor3}. We recall some interesting results, which will be used in our paper.

\begin{defn}\cite{Jor3}\label{exreg}{\rm
 For any $M \in \mathrm{D}(A)$, we define the Ext-regularity of $M$ by
$$\mathrm{Ext.reg}M = -\inf\{i|H^i(R\Hom_{A}(M,\k))\neq 0\},$$ and similarly for $M \in \mathrm{D}(A^{op})$.  Note that $\mathrm{Ext.reg}(0)=-\infty$.  }
\end{defn}

Let $A$ be a homologically smooth connected cochain DG algebra.  Then both ${}_A\k$ and $\k_A$ are compact.  Let $K$ and $L$ be minimal semi-free resolutions of ${}_A\k$ and $\k_A$, respectively.  We have $\langle K\rangle = \langle {}_A\k \rangle$ and $\langle L\rangle =\langle \k_A\rangle$  in $\mathrm{D}(A)$. In this case, we have \cite[Setup 4.1]{Jor3}.

Set $\mathcal{N}=\langle {}_A\k\rangle^{\bot} =\langle {}_AK\rangle^{\bot}$ in $\mathrm{D}(A)$, and define the torsion and the complete DG modules by $\mathrm{D}^{\mathrm{tors}}(A)={}^{\bot}\mathcal{N}, \mathrm{D}^{\mathrm{comp}}(A)=\mathcal{N}^{\bot}$. Then
 $$\mathrm{D}^{\mathrm{tors}}(A)=\langle{}_A\k\rangle=\langle{}_AK\rangle. $$
Let $\mathcal{E}=\Hom_{A}(K,K)$ be the endomorphism DG algebra. The DG module $K$ acquires the structure ${}_{A,{}_{\mathcal{E}}}K$ while $K^*=\Hom_{A}(K,A)$ has the structure $K^*_{A,\mathcal{E}}$. Define functors
\begin{align*}
T(-)=-\,{}^L{\otimes}_{\mathcal{E}}K ,\\
R(-)=\Hom_A(K,-)\simeq K^*\,{}^L{\otimes}_A- ,\\
C(-)=R\Hom_{\mathcal{E}^{op}}(K^*,-) ,
\end{align*}
which form adjoint pairs $(T,R)$ and $(R,C)$ between $\mathrm{D}(\mathcal{E}^{op})$ and $\mathrm{D}(A)$.
There are pairs of quasi-inverse equivalences of categories as follows
$$\xymatrix{\mathrm{D}^{\mathrm{comp}}(A)   \ar@<1ex>[r]^{\quad R}
& \mathrm{D}(\mathcal{E}^{op}) \ar@<1ex>[l]^{\quad C}  \ar@<1ex>[r]^{T} & \mathrm{D}^{\mathrm{tors}}(A) \ar@<1ex>[l]^{ R} } .$$
In particular, $RC$ and $RT$ are equivalent to the identity functor on $\mathrm{D}(\mathcal{E}^{op})$ , so if we set $$\Gamma =TR, \Lambda=CR,$$
then we get endofunctors of $\mathrm{D}(A)$ which form an adjoint pair $(\Gamma, \Lambda)$ and satisfy
$$\Gamma^2\simeq \Gamma, \Lambda^2\simeq \Lambda, \Gamma\Lambda\simeq \Gamma, \Lambda\Gamma\simeq \Lambda.$$
These functors are adjoints of inclusions as follows, where left-adjoint are displayed above right-adjoints
$$\xymatrix{\mathrm{D}^{\mathrm{comp}}(A)   \ar@<1ex>[r]^{\quad \mathrm{inc}}
& \mathrm{D}(A) \ar@<1ex>[l]^{\quad \Lambda}  \ar@<1ex>[r]^{\Gamma} & \mathrm{D}^{\mathrm{tors}}(A) \ar@<1ex>[l]^{ \mathrm{inc }} } .$$
Write $F=K^*{}^L\otimes_{\mathcal{E}}K$ and $D=F^{\vee}=\Hom_{\k}(F,\k)$. One sees that $F$ and $D$ have the structures ${}_AF_A$ and ${}_AD_A$,  respectively.  From the definitions, we have $$\Gamma(-)=F\,{}^L{\otimes}_A-\quad \text{and} \quad \Lambda(-)=R\Hom_A(F,-).$$
\begin{defn}\cite[Definition 5.1]{Jor3}{\rm
For $M\in \mathrm{D}(A)$, we define the Castelnuovo-Mumford regularity of $M$ by
$$ \mathrm{CMreg}M=\sup\{i|H^i(\Gamma (M))\neq 0\}.$$
Note that $ \mathrm{CMreg}(0)=-\infty$. }
\end{defn}

\begin{lem}\cite[Proposition 5.6]{Jor3}\label{cmformula}{\rm  Let $A$ be a homologically smooth connected cochain DG algebra. If $M$ is an object in $\mathrm{D}^+(A)$. Then
\begin{enumerate}
\item $\mathrm{CMreg}M\neq -\infty$.
\item $\mathrm{Ext.reg}M\le \mathrm{CMreg}M+\mathrm{Ext.reg}\k$.
\item $\mathrm{CMreg}M\le \mathrm{Ext.reg}M+\mathrm{CMreg}A$.
\end{enumerate}
}
\end{lem}

\begin{lem}\cite[(5.3)]{Jor3}\label{cmfin}
{\rm Let $A$ be a homologically smooth connected cochain DG algebra such that $H(A)$ is Noetherian with a balanced dualizing complex. Then for any non-zero object in $\mathrm{D}^f(A)$, we have $-\infty< \mathrm{CMreg}M<+\infty$. }
\end{lem}

\begin{lem}\cite[Theorem 5.7]{Jor3}\label{mkoszul}{\rm
Assume that $A$ is a homologically smooth connected cochain DG algebra such that $H(A)$ is Noetherian with a balanced dualizing complex. Let $M$ be a non-zero object in $\mathrm{D_{fg}}(A)$.
\begin{enumerate}
\item If $\mathrm{Ext.reg}\k<\infty$, then $\mathrm{Ext.reg}M<\infty$.

\item If $A$ is a Koszul DG algebra and $\mathrm{CMreg}M\le t$ for some $t$, then $M^{\ge t}$ is a Koszul DG module.
\end{enumerate}}
\end{lem}

%%% ----------------------------------------------------------------------
\section{dg free class of semi-free DG modules}
The terminology `class'
in group theory is used to measure the shortest length of a
filtration with sub-quotients of certain type. Carlsson \cite{Car}
introduced `free class' for solvable free differential graded
modules over a graded polynomial ring. In \cite{ABI},
 Avramov, Buchweitz and Iyengar  introduced free class,
projective class and flat class for differential modules over a
commutative ring. Inspired from them, the author and Wu \cite{MW3} introduced the
the DG free class for
semi-free DG modules.

\begin{defn}\cite{MW3}{\rm Let $F$ be a semi-free DG
$A$-module. A semi-free filtration of $F$
$$0=F(-1)\subseteq F(0)\subseteq\cdots \subseteq F(n)\subseteq\cdots$$
is called strictly increasing, if $F(i-1) \neq F(i)$ when
$F(i-1)\neq F, i\ge 0$.  If there is some $n$ such that $F(n)=F$ and
$F(n-1)\neq F$, then we say that this strictly increasing semi-free
filtration has length $n$. If no such integer exists, then we say
the length is $+\infty$.}
\end{defn}

\begin{defn}\cite{MW3}{\rm
Let $F$ be a semi-free DG $A$-module. The DG free class of $F$ is
defined to be the number
$$ \inf\{n\in \Bbb{N}\cup \{0\}\,|\, F \,\text{admits a strictly
increasing semi-free filtration of length}\,\, n\}.$$ We denote it
as $\mathrm{DG free\,\, class}_AF$. }
\end{defn}

\begin{defn}\cite{MW3}{\rm
Let $M$ be a DG $A$-module. The cone length of $M$ is defined to be
the number
$$\mathrm{cl}_AM =
\inf\{\mathrm{\,DGfree\,\,class}_AF\,|\,F \stackrel{\simeq}\to M
 \ \text{is a semi-free resolution of}\  M\}.$$ }
\end{defn}
Cone length of a DG $A$-module plays a similar role in DG
homological algebra as projective dimension of a module over a ring
does in classic homological ring theory (cf. \cite{MW3}).
\begin{rem}\label{conefin} Note that $\mathrm{cl}_AM$ may be $+\infty$.  We call this invariant
`cone length' because semi-free DG $A$-modules can be constructed by
iterative cone constructions from DG free $A$-modules (see \cite[Lemma
4.1]{MW3}) and any DG $A$-module admits a semi-free resolution.
For any compact DG $A$-module $M$, it admits a minimal semi-free resolution $F_M$, which has a finite semi-basis.
We have $\mathrm{\,DGfree\,\,class}_AF_M<\infty$ and then  $\mathrm{cl}_AM<\infty$.
\end{rem}

By \cite[Proposition 2.4]{MW1}, any object $M\in \mathrm{D_{fg}}(A)$ admits a minimal semi-free resolution $F_M$. A natural question is whether $\mathrm{cl}_AM=\mathrm{DGfree.class}_AF_M$.  One will see that this question is closely related to Question \ref{mainquest} from the proof of Theorem \ref{koszulcomp}.  Actually, if $$\mathrm{cl}_AM=\mathrm{DGfree.class}_AF_M,  \forall  M\in \mathrm{D_{fg}}(A),$$ then one has $\mathrm{D_{fg}}(A)=\mathrm{D}^c(A)$, when $A$ is a homologically smooth connected cochain DG algebra and $H(A)$ is Noetherian.
By \cite[Proposition 3.6]{MW3}, we have $\mathrm{cl}_AM=\mathrm{DGfree.class}_AF_M$ when $M$ admits a minimal Eilenberg-Moore resolution. However, we are unable to prove this in general. In spite of this, we do make some progress concerning Question \ref{mainquest} by Theorem A and Theorem B. In order to prove these two theorems, we need to make the following preparations.

\begin{lem}\cite[\S 12, Theorem 3.2]{AFH}\label{semidecomp}
 Let $F$ be a semi-free DG $A$-module such that $H(F)$ is bounded
below. Then there is a minimal semi-free resolution $G$ of $F$ and a
homotopically trivial DG $A$-module $Q$ such that $F \cong G\oplus
Q$ as a DG $A$-module.
\end{lem}

\begin{lem}\cite[\S 8, Corollary 4.3]{AFH}\label{semihomotopy}
 Let $f: F \to F'$ be a morphism between two semi-free DG $A$-modules. Then $f$ is a quasi-isomorphism if and only if it is a homotopy equivalence.
\end{lem}

\begin{prop}\label{finitecl}
Let $A$ be a connected DG algebra such that
$\mathrm{cl}_{A^e}A< +\infty$. Then for any DG $A$-module $M$,  we
have $\mathrm{cl}_AM\le \mathrm{cl}_{A^e}A<+\infty$. Morever, if $A$ is homologically smooth and $M$ is an object in $\mathrm{D}^+(A)$, then $M$ admits a bounded below semi-free resolution $F$ such that $\inf\{i|F^i\neq 0\}=\inf\{i|H^i(M)\neq 0\}$ and $$\mathrm{DGfree\,\,class}_AF\le \mathrm{DGfree\,\,class}_{A^e}P,$$ where $P$ is a minimal semi-free resolution of ${}_{A^e}A$.
\end{prop}
\begin{proof}
Let $\mathrm{cl}_{A^e}A=n$. By the definition of cone length, the DG
$A^e$-module $A$ admits a semi-free resolution $X$ such that
$\mathrm{DGfree\,\,class}_{A^e}X = n$. This implies that $X$ admits
a strictly increasing semi-free filtration
$$0\subset X(0)\subset X(1)\subset\cdots \subset X(n)=X,$$
where $X(0)= A^e\otimes V(0)$ and each $X(i)/X(i-1)\cong A^e\otimes
V(i)$ is a DG free $A^e$-module. Let $E_i=\{e_{i_j}|j\in
I_i\}, i \ge 0,$ be a basis of $V(i)$. For any $i\ge 1$, define
$f_i: A^e\otimes \Sigma^{-1}V(i)\to X(i-1)$ such that
$f_i(\Sigma^{-1}e_{i_j}) =
\partial_{X(i)}(e_{i_j})$. By \cite[Lemma 4.1]{MW3}, $X(i) \cong
\mathrm{cone}(f_i), i=1, 2,\cdots, n$.

For any DG $A$-module $M$, let $\varrho_M:F_M\to M$ be a semi-free resolution of $M$. As
a DG $A$-module, $X(i)\otimes_AF_M \cong \mathrm{cone}(f_i\otimes_A
\mathrm{id}_ {F_M})$, $i=1, 2,\cdots, n$.
 Since $A^e\otimes_AF_M\cong A\otimes F_M$, we have
$$(A^e\otimes V(i))\otimes_AF_M\cong A\otimes V(i)\otimes F_M,
\quad i= 0,1,\cdots, n.$$ Choose a subset $\{m \}\subseteq F_M$ such
that each $m$ is a cocycle and $\{\lceil m\rceil \}$ is a basis of the
$\k$-vector space $H(F_M)$. Define a DG morphism
$$\phi_i: A\otimes V(i)\otimes H(F_M)\to
A\otimes V(i)\otimes F_M$$ such that $\phi_i(a \otimes v \otimes
\lceil m\rceil)=a \otimes v \otimes m$, for any $a \in A, v \in V(i)$ and
$\lceil m\rceil$. It is easy to check that $\phi_i$ is a quasi-isomorphism.

In the following, we prove inductively that $\mathrm{cl}_A
(X(i)\otimes_AF_M)\le i, i= 0,1,\cdots, n$. Since $\phi_0: A\otimes
V(0)\otimes H(F_M)\to X(0)\otimes_AF_M$ is a quasi-isomorphism, we
have $\mathrm{cl}_A (X(0)\otimes_A F_M) = 0$. Suppose inductively
that we have proved that $$\mathrm{cl}_A(X(l)\otimes_AF_M)\le l,\,
l\ge 0.$$ We should prove $\mathrm{cl}_A(X(l+1)\otimes_AF_M)\le
l+1$. Since $\mathrm{cl}_A(X(l)\otimes_AF_M)\le l$, there is a
semi-free resolution $\varphi_l:F_l\stackrel{\simeq}{\to}
X(l)\otimes_AF_M$ such that $\mathrm{DGfree\,\,class}_AF_l\le l$.
 Because $A\otimes\Sigma^{-1}V(l+1)\otimes H(F_M)$ is semi-free,
there is a DG morphism
$$\psi_{l}:A\otimes\Sigma^{-1}V(l+1)\otimes H(F_M)\to
F_l$$ such that $\varphi_l\circ \psi_l \sim
(f_l\otimes_A\mathrm{id}_{F_M})\circ \Sigma^{-1}(\phi_{l+1})$.

 For convenience, we write
$Q(l+1)=A\otimes V(l+1)$ and $K(l+1) =A^e\otimes V(l+1)$. In
$\mathrm{D}(A)$, there is a morphism
$h_{l+1}:\mathrm{cone}(\psi_l)\to X(l+1)\otimes_AF_M$ making the
diagram
\begin{tiny}
\begin{align*}
\xymatrix{\Sigma^{-1}Q(l+1)\otimes H(F_M)
\ar[r]^{\psi_l}\ar[d]^{\Sigma^{-1}(\phi_{l+1})} & F_l
\ar[d]^{\varphi_l}\ar[r]^{\tau_{l}}&\mathrm{cone}(\psi_l)\ar[d]^{\exists
h_{l+1}}\ar[r]^{\varepsilon_l}&Q(l+1)\otimes_{\k}H(F_M)\ar[d]^{\phi_{l+1}}
\\
\Sigma^{-1}K(l+1)\otimes_AF_M \ar[r]^{f_l\otimes_A\mathrm{id}_{F_M}}
&X(l)\otimes_AF_M \ar[r]^{\iota_l}&X(l+1)\otimes_AF_M
\ar[r]^{\pi_l}& K(l+1)\otimes_AF_M
\\}
\end{align*}
\end{tiny}
commute. By five-lemma, $h_{l+1}$ is an isomorphism in
$\mathrm{D}(A)$. This implies that there are quasi-isomorphisms
$g:Y\to \mathrm{cone}(\psi_l)$ and $t:Y\to X(l+1)\otimes_AF_M$,
where $Y$ is some DG $A$-module. Hence $ \mathrm{cl}_A
(X(l+1)\otimes_AF_M)=\mathrm{cl}_AY =
\mathrm{cl}_A\mathrm{cone}(\psi_l)\le l+1 $.

We have proved inductively that $\mathrm{cl}_A(X\otimes_AF_M)\le n$.
Since $F_M\simeq X\otimes_AF_M$, we get $\mathrm{cl}_AM\le n$.

Now, assume that $A$ is homologically smooth and $M$ is an object in $\mathrm{D}^+(A)$.
  Let $\alpha:P \to {}_{A^e}A$ be a minimal semi-free resolution and $b=\inf\{i|H^i(M)\neq 0\}$. By \cite[Proposition 2.4]{MW1}, we have $\inf\{i|P^i\neq 0\}=0$ and we may assume that $\varrho_M$ mentioned above is a minimal semi-free resolution of $M$.
Since ${}_{A^e}A$ is compact, $P$ admits a finite semi-basis. So $\mathrm{DGfree\,\,class}_{A^e}P$ is finite. Let $ \mathrm{DGfree\,\,class}_{A^e}P=t$. Then $P$ admits
a strictly increasing semi-free filtration
$$0\subset P(0)\subset P(1)\subset\cdots \subset P(t)=P,$$
where $P(0)= A^e\otimes W(0)$ and each $P(i)/P(i-1)\cong A^e\otimes
W(i), i\ge 1$, is a DG free $A^e$-module. Let $\Lambda_i=\{\lambda_{i_j}|j\in
J_i\}, i \ge 0,$ be a basis of $W(i)$. For any $i\ge 1$, define
$g_i: A^e\otimes_{\k} \Sigma^{-1}W(i)\to P(i-1)$ such that
$g_i(\Sigma^{-1}\lambda_{i_j}) =
\partial_{P(i)}(\lambda_{i_j})$. By  \cite[Lemma 4.1]{MW3}, $P(i) \cong
\mathrm{cone}(g_i), i=1, 2,\cdots, t$.
As
a DG $A$-module, $P(i)\otimes_AF_M \cong \mathrm{cone}(g_i\otimes_A
\mathrm{id}_ {F_M})$, $i=1, 2,\cdots, t$. Define a DG morphism
$$\eta_i: A\otimes W(i)\otimes H(F_M)\to
A\otimes W(i)\otimes F_M$$ such that $\eta_i(a \otimes v \otimes
\lceil m\rceil )=a \otimes v \otimes m$, for any $a \in A, v \in V(i)$ and
$\lceil m\rceil$. It is easy to check that $\eta_i$ is a quasi-isomorphism.
Since $P(0)\otimes_AF_M=A\otimes W(i)\otimes F_M$ and $\eta_0:A\otimes W(0)\otimes H(F_M)\to
A\otimes W(0)\otimes F_M$ is a semi-free resolution with $$\mathrm{DGfree\,\,class}_AA\otimes W(0)\otimes H(F_M)=0.$$ Let $G_0= A\otimes W(0)\otimes H(F_M)$. Then $\inf\{j|G_0^i\neq 0\}=b$.  We assume inductively that $P(i)\otimes_AF_M $ admits a semi-free resolution $\omega_i:G_i\to P(i)\otimes_AF_M $ with
$b=\inf\{j|G_i^j\neq 0\}$ and $\mathrm{DGfree\,\,class}_{A}G_i \le i$ when $i\le r <t$. We want to show that
 $P(r+1)\otimes_AF_M$ admits a semi-free resolution $G_{r+1}$ with
$b=\inf\{j|G_{r+1}^j\neq 0\}$ and $\mathrm{DGfree\,\,class}_{A}G_{r+1} \le r+1$.
 Since $A\otimes \Sigma^{-1}W(r+1)\otimes H(F_M)$ is semi-free,
there is a DG morphism
$$\beta_{r}:A\otimes \Sigma^{-1}W(r+1)\otimes H(F_M)\to
G_r$$ such that $\omega_r\circ \beta_r \sim
(g_r\otimes_A\mathrm{id}_{F_M})\circ\Sigma^{-1}(\eta_{r+1})$.
Let $R(l+1)=A\otimes W(l+1)$ and $T(l+1) =A^e\otimes W(l+1)$. In
$\mathrm{D}(A)$, there is a morphism
$$\theta_{r+1}:\mathrm{cone}(\beta_r)\to P(r+1)\otimes_AF_M$$ making the
diagram
\begin{tiny}
\begin{align*}
\xymatrix{\Sigma^{-1}R(r+1)\otimes H(F_M)
\ar[r]^{\beta_r}\ar[d]^{\Sigma^{-1}(\eta_{r+1})} & G_r
\ar[d]^{\omega_r}\ar[r]^{\nu_{r}}&\mathrm{cone}(\beta_r)\ar[d]^{\exists
\theta_{r+1}}\ar[r]^{\gamma_r}&R(r+1)\otimes H(F_M)\ar[d]^{\eta_{r+1}}
\\
\Sigma^{-1}T(r+1)\otimes_AF_M \ar[r]^{g_r\otimes_A\mathrm{id}_{F_M}}
&P(r)\otimes_AF_M \ar[r]^{\iota_r}&P(r+1)\otimes_AF_M
\ar[r]^{\pi_r}& T(r+1)\otimes_AF_M
\\}
\end{align*}
\end{tiny}
\!commute. By five-lemma, $\theta_{r+1}$ is an isomorphism in
$\mathrm{D}(A)$. This implies that there are quasi-isomorphisms
$\rho:Z\to \mathrm{cone}(\beta_r)$ and $\chi:Z\to P(r+1)\otimes_AF_M$,
where $Z$ is some DG $A$-module. It is easy to see that $\inf\{i|\mathrm{cone}(\beta_r)^i\neq 0\}=b$.
So $Z\in \mathrm{D}^+(A)$ and $Z$ has a minimal semi-free resolution $\xi:F_Z\to Z$ by \cite[Proposition 2.6]{MW1}.
Since $\rho\circ \xi:  F_Z\to \mathrm{cone}(\beta_r)$ is a quasi-isomorphism between semi-free DG modules,
 $\rho\circ \xi$ is a homotopy equivalence by Lemma \ref{semihomotopy}. Let $\sigma: \mathrm{cone}(\beta_r)\to F_Z$ be its homotopy inverse. Then $$\chi\circ \xi \circ \sigma: \mathrm{cone}(\beta_r)\to P(r+1)\otimes_AF_M$$ is a quasi-isomorphism and it is a semi-free resolution of  $P(r+1)\otimes_AF_M$. Let $G_{r+1}=\mathrm{cone}(\beta_r)$. Since $\mathrm{DGfree\,\,class}_{A}G_{r}\le r$ by the induction hypothesis, we have $\mathrm{DGfree\,\,class}_{A}G_{r+1}\le r+1$.
By the induction above, $P(t)\otimes_AF_M$ admits a semi-free resolution $$\theta_t:F=G_{t}\to P(t)\otimes_AF_M $$ such that
$\inf\{i|F^i\neq 0\}=b$ and $\mathrm{DGfree\,\,class}_{A}F\le t$. Since the composition
$$P(t)\otimes_AF_M\stackrel{\alpha\otimes \mathrm{id}_{F_M}}{\to} {}_{A^e}A\otimes_AF_M\stackrel{\varrho_M}{\to}M$$
is a quasi-isomorphism, $F$ is a semi-free resolution of $M$.
\end{proof}

\begin{prop}\label{kosmin}
Let $A$ be a homologically smooth connected cochain DG algebra. Then any Koszul DG module $M$ admits a minimal semi-free resolution $F_M$ with $$\mathrm{DGfree\,\,class}_{A}F_M\le \mathrm{DGfree\,\,class}_{A^e}P<\infty,$$ where $P$ is a minimal semi-free resolution of ${}_{A^e}A$.
\end{prop}
\begin{proof}
By Proposition \ref{finitecl},
$M$ has a bounded below semi-free resolution $F$ such that $\inf\{i|F^i\neq 0\}=\inf\{i|H^i(M)\neq 0\}=b$ and $$\mathrm{DGfree\,\,class}_AF\le \mathrm{DGfree\,\,class}_{A^e}P<\infty,$$ where $P$ is a minimal semi-free resolution of ${}_{A^e}A$. By the proof of Proposition \ref{finitecl}, one sees that $F$ is not minimal in general.
Since $M$ is Koszul, it admits a minimal semi-free resolution $F_M$, which has a semi-basis concentrated in degree $b$ by Remark \ref{onkoszul}. We have $F\cong F_M\oplus Q$ by Lemma \ref{semidecomp}, where $Q$ is a homotopically trivial DG $A$-module. Let $L=F_M\oplus Q$ and  $\mathrm{DGfree\,\,class}_AF=t$. Then $\mathrm{DGfree\,\, class}_AL =t$ and $\inf\{i|L^i\neq 0\}=b$. So $L$ admits a semi-free filtration
$$0=L(-1)\subset L(0)\subset L(1)\subset \cdots \subset L(t)=L.$$ Then $L$ admits a semi-basis $E=\bigsqcup\limits_{i=0}^t E_i$ such that $\partial_L(E_i)\subseteq A(\bigsqcup\limits_{j=0}^{i-1}E_j) $ for any $i\in \{1,2,\cdots, t\}$. Let $$E_i=\{e_{\lambda}|\lambda\in I_i \},\quad V_i=\bigoplus_{e\in E_i} ke
\quad \text{and}\quad V=\bigoplus\limits_{i=0}^t V_i.$$
Then $L^{\#}=A^{\#}\otimes V$.
Since $L^{\#}=F_M^{\#}\oplus Q^{\#}$ is a bounded below free graded $A^{\#}$-module, $Q^{\#}$ is a bounded below projective graded $A^{\#}$-module. So $Q^{\#}$ is a free graded $A^{\#}$-module. We may assume that $F_M^{\#}=A^{\#}\otimes W$ and $Q^{\#}=A^{\#}\otimes U$. Then $A^{\#}\otimes V = (A^{\#}\otimes W )\oplus (A^{\#}\otimes U)$. We have $V^b=W\oplus U^b$ by considering the degree. Let $V^b_i=V^b\bigcap V_i, i=0,1,\cdots, t$. Then $\bigcup\limits_{i=0}^tV^b_i=\bigcup\limits_{i=0}^t(V^b\bigcap V_i)= V^b\bigcap (\bigcup\limits_{i=0}^t V_i)=V^b$.
Let $I^b_i\subset I_i$ such that $E^b_i=\{e_{\lambda}|\lambda\in I^b_i\}$ is a basis of $V^b_i$. For any $\lambda\in I^b_i$, there are $\omega_{\lambda}\in W$ and $u_{\lambda}\in U^b$ such that $e_{\lambda}=\omega_{\lambda}+u_{\lambda}$. Then for any $\lambda\in I^b_i$, we have $\partial_L(e_{\lambda})\in L(i-1)^{b+1}$ and hence
\begin{align*}
\partial_{F_M}(\omega_{\lambda})+\partial_Q(u_{\lambda})=\partial_L(e_{\lambda})&=\sum\limits_{j=0}^{i-1}\sum\limits_{\lambda\in I^b_j}a_{\lambda}e_{\lambda}+\chi_{\lambda}\\
&=\sum\limits_{j=0}^{i-1}\sum\limits_{\lambda\in I^b_j}a_{\lambda}\omega_{\lambda}+\sum\limits_{j=0}^{i-1}\sum\limits_{\lambda\in I^b_j}a_{\lambda}u_{\lambda} +\chi_{\lambda},
\end{align*}
where $a_{\lambda}\in \mathcal{A}^1$ and $\chi_{\lambda}\in (\mathcal{A}^{\#}\otimes V^{\ge b+1})\bigcap L(i-1)^{b+1}$.
One sees  that $\chi_{\lambda}\in Q$ since $W$ is concentrated in degree $b$. Hence
\begin{align}\label{diffcond}
\partial_{F_M}(\omega_{\lambda})=\sum\limits_{j=0}^{i-1}\sum\limits_{\lambda\in I^b_j}a_{\lambda}\omega_{\lambda}
\end{align}
and
$$\partial_Q(u_{\lambda})= \sum\limits_{j=0}^{i-1}\sum\limits_{\lambda\in I^b_j}a_{\lambda}u_{\lambda} +\chi_{\lambda}.$$
Let $W_{\le j}=\sum\limits_{i=0}^j\sum\limits_{\lambda\in I^b_i}\k\omega_{\lambda}, j=0,1,\cdots t$. By (\ref{diffcond}),
each $A\otimes W_{\le j}$ is a DG submodule of $F_M$ and
\begin{align*}
0\subseteq A\otimes W_{0}\subseteq A\otimes W_{\le 1}\subseteq \cdots \subseteq A\otimes W_{\le t}=F_M
\end{align*}
is a semi-free filtration of $F_M$. Therefore, $\mathrm{DGfree.class}_AF_M\le t$.

\end{proof}

\section{two main theorems}
 In this section, we present the main results of this paper.
 Applying Proposition \ref{kosmin}, we
  can prove the following theorem, which partially solves Question \ref{mainquest}.

\begin{thm}\label{koszulcomp}
Let $A$ be a homologically smooth connected cochain DG algebra such that $H(A)$ is Noetherian. Then
any Koszul DG module in $\mathrm{D_{fg}}(A)$ is compact.
\end{thm}

\begin{proof}
By Proposition \ref{kosmin}, $M$ admits a minimal semi-free resolution $F_M$ such that  $$\mathrm{DGfree\,\,class}_{A}F_M\le \mathrm{DGfree\,\,class}_{A^e}P<\infty,$$
where $P$ is a minimal semi-free resolution of ${}_{A^e}A$. Let $\mathrm{DGfree\,\,class}_{A}F_M=t$. Then there is a semi-free filtration
$$0=F_M(-1)\subset F_M(0)\subset F_M(1)\subset \cdots \subset F_M(t)=F_M$$
of $F_M$. Each $F_M(i)/F_M(i-1) = A\otimes W_i$ is
a DG free $A$-module on a cocycle basis, $i=0,1,\cdots, t$. We should prove that each
$W_i$ is finite dimensional. Let $\{e_{i,j}|j\in I_i\}$ be a basis of
$W_i, i=0,1, \cdots, t$. Let $\iota_0: F_M(0)\to F_M$ be the
inclusion map. Since $\mathrm{im}H(\iota_0)$ is a graded
$H(A)$-submodule of $H(F_M)$ and $H(A)$ is Noetherian,
$\mathrm{im}H(\iota_0)\cong
\frac{H(F_M(0))}{\mathrm{ker}H(\iota_0)}$ is also a finitely
generated $H(A)$-module.  Let $\mathrm{im}H(\iota_0) =
H(A)f_{0,1}+H(A)f_{0,2}+\cdots+H(A)f_{0,n}$. Since $H(F_M(0))\cong
\bigoplus\limits_{j\in I_0}H(A)e_{0,j}$ is a free graded
$H(A)$-module, there is a finite subset $J_0=\{i_1,i_2,\cdots,i_l\}$
of $I_0$ such that
$$f_{0,s}=\sum\limits_{t=1}^la_{st}\overline{e_{0,i_t}}, s=1,2,
\cdots, n,$$ where each $a_{st}\in H(A)$. If $V(0)$ is infinite
dimensional, then both $I_0$ and $I_0\setminus J_0$ are infinite
sets. Hence for any $j\in I_0\setminus J_0$, we have $e_{0,j}\in
\mathrm{ker}H(\iota_0)$. Since $[\iota_{0}(e_{0,j})]=[e_{0,j}]=0$ in
$H(F_M)$, there exists $x_{0,j}\in F_M$ such that
$\partial_{F_M}(x_{0,j})=e_{0,j}$. This contradicts with the
minimality of $F_M$. Thus $W_0$ is finite dimensional and
$F_M(0)\in \mathrm{D_{fg}}(A)$.

Suppose inductively that we have proved that $W_0, \cdots, W_{i-1}$
are finite dimensional. Then each $F_M(j)/F_M(j-1), (j=0,1,\cdots,
i-1)$ is an object in $\mathrm{D_{fg}}(A)$. And we can prove inductively that
each $F_M(j), (j=0,1,\cdots, i-1)$ is in $\mathrm{D_{fg}}(A)$ by the
following sequence of short exact sequences
\begin{align*}
0\longrightarrow F_M(j-1)\longrightarrow F_M(j)\longrightarrow
F_M(j)/F_M(j-1)\longrightarrow 0 , j= 1,\cdots, i-1.\\
\end{align*}
Similarly, $F_M/F_M(i-1)$ is also an object in $\mathrm{D_{fg}}(A)$ by
the short exact sequence $$0\longrightarrow F_M(i-1) \longrightarrow
F_M \longrightarrow F_M/F_M(i-1)\longrightarrow 0.$$ On the other
hand, it is easy to see that $F_M/F_M(i-1)$ is also a minimal
semi-free DG $A$-module and it has a semi-free filtration
$$F_M(i)/F_M(i-1)\subseteq F_M(i+1)/F_M(i-1)\subseteq \cdots \subseteq
F_M(n)/F_M(i-1)=F_M/F_M(i-1).$$ Let $\iota_i: F_M(i)/F_M(i-1)\to
F_M/F_M(i-1)$ be the inclusion morphism. Since
$\mathrm{im}H(\iota_i)$ is a graded $H(A)$-submodule of
$H(F_M/F_M(i-1))$ and $H(A)$ is Noetherian,
$\mathrm{im}H(\iota_i)\cong
\frac{H(F_M(i)/F_M(i-1))}{\mathrm{ker}H(\iota_i)}$ is also a
finitely generated $H(A)$-module.  Let $\mathrm{im}H(\iota_i) =
H(A)f_{i,1}+H(A)f_{i,2}+\cdots+H(A)f_{i,m}$. Since
$$H(F_M(i)/F_M(i-1))\cong \bigoplus\limits_{j\in I_i}H(A)e_{i,j}$$ is a
free graded $H(A)$-module, there is a finite subset
$J_i=\{s_1,s_2,\cdots,s_r\}$ of $I_i$ such that
$$f_{i,l}=\sum\limits_{t=1}^ra_{l,t}\overline{e_{i,s_t}}, l=1,2,
\cdots, m,$$ where each $a_{lt}\in H(A)$. If $W_i$ is an infinite
dimensional space, then both $I_i$ and $I_i\setminus J_i$ are
infinite sets. Hence for any $j\in I_i\setminus J_i$, we have
$e_{i,j}\in \mathrm{ker}H(\iota_i)$. Since
$[\iota_{i}(e_{i,j})]=[e_{i,j}]=0$ in $H(F_M/F_M(i-1))$, there exist
$x_{i,j}\in F_M/F_M(i-1)$ such that
$\partial_{F_M}(x_{i,j})=e_{i,j}$. This contradict with the
minimality of $F_M$. Thus $W_i$ is finite dimensional.

By the induction above, we prove that each $W_i, (i=0,1,\cdots,n)$
is finite dimensional. Hence $F_M$ has a finite semi-basis and $M$
is compact.

\end{proof}
Proposition \ref{koszulcomp} partially solves Question \ref{mainquest}, since we add Koszul condition on DG modules in $\mathrm{D_{fg}}(A)$.
Inspired from \cite[Theorem 5.7]{Jor3}, we pass on the Koszul condition to homologically smooth DG algebras. We want to show that
   $\mathrm{D_{fg}}(A)=\mathrm{D}^c(A)$ when $A$ is a Koszul homologically smooth DG algebra such that $H(A)$ is Noetherian and admits a balanced dualizing complex. We need first to prove the following lemma and propostion.

\begin{lem}\label{quasitrivial} {\rm  Let $A$ be a connected cochain DG algebra. Assume that $F_{\k}$ is a minimal semi-free resolution of $\k_{A}$.
If $M$ is a bounded below DG $A$-module such that $F_{\k}\otimes_A M$ is quasi-trivial, then $M$ is quasi-trivial. }
\end{lem}
\begin{proof}
If $M$ is not quasi-trivial, then there is $b\in \Bbb{Z}$ such that  $b=\inf\{i|H^i(M)\neq 0\}$.
By \cite[Proposition 2.4]{MW1}, $M$ admits a minimal semi-free resolution $F_M$ with $F^{\#}_M\cong
\coprod_{i\ge b}\Sigma^{-i}(A^{\#})^{(\Lambda^i)}$
each $\Lambda^i$ is an index set. Then \begin{align*}
H(F_{\k}\otimes_A M)=H(\k\,^L{\otimes}_AM)=H(\k\otimes_AF_M)= \coprod_{i\ge b}\Sigma^{-i}(\k)^{(\Lambda^i)}\neq 0.
\end{align*}
This contradicts with the condition that $F_{\k}\otimes_A M$ is quasi-trivial. Hence $M$ is quasi-trivial.
\end{proof}

\begin{prop}\label{fgcomp}
Let $A$ be a Koszul and homologically smooth connected  cochain DG algebra such that $H(A)$ is
a Noetherian graded algebra with a balanced dualizing complex. Then any non-zero object in $\mathrm{D_{fg}}(A)$ is compact.

\end{prop}
\begin{proof}
Since $A$ is a non-zero object in $\mathrm{D_{fg}}(A)$, we have $-\infty< \mathrm{CMreg}A<\infty$
by Lemma \ref{cmfin}.
For any non-zero object $M\in \mathrm{D_{fg}}(A)$, we have $$\mathrm{CMreg}M\le \mathrm{Ext.reg}M+\mathrm{CMreg}A$$ by Lemma \ref{cmformula}.
By the assumption that $A$ is Koszul, we have $\mathrm{Ext.reg}_A\k=0<\infty$ and then $\mathrm{Ext.reg}M<\infty$ by Lemma \ref{mkoszul}. Hence $\mathrm{CMreg}M $ is finite and $\mathrm{CMreg}M\le t$ for some $t$. By Lemma \ref{mkoszul}, $M^{\ge t}$ is a Koszul DG module.

For any $i\in \Bbb{Z}$, the $\k$-vector space $M^{i}$
can be decomposed as $B^i \oplus H^i \oplus C^i$ with $B^i\oplus
H^i = \mathrm {Ker}\, \partial^i_M, C^i\stackrel{\cong}{\to} B^{i+1} =
\mathrm {Im} \, \partial^i_M$ and $H^i \cong H^i(M)$. We get a DG $A$-submodule
\begin{align*}
\tau^{\ge i}M:\qquad\qquad  \cdots\to 0\to H^i\oplus C^i\stackrel{\partial_M^i}{\to} M^{i+1}\stackrel{\partial_M^{i+1}}{\to}
M^{i+2}\stackrel{\partial_M^{i+2}}{\to}\cdots
\end{align*}
of $M$. We have $$H^j(\tau^{\ge i}M)=\begin{cases}
0,   \quad\quad \quad\quad  j<i \\
H^j(M), \quad \,\, j\ge i.
\end{cases}$$
Hence $H(\tau^{\ge u}M)$ is a graded $H(A)$-submodule of $H(M)$. Since $M\in \mathrm{D_{fg}}(A)$ and $H(A)$ is a Noetherian graded algebra, $H(\tau^{\ge u}M)$ is a finitely generated graded $H(A)$-module.  So $\tau^{\ge i}M \in \mathrm{D_{fg}}(A)$.
Let $b=\inf\{i|H^i(M)\neq 0\}$. We claim that there exist some $i\ge b$ such that $\tau^{\ge i}M$ is a quasi-trivial DG $A$-module or a compact DG $A$-module.

If $H^t(M^{\ge t})=0$, then there exists $u > t$  such that $M^{\ge t}$ has a minimal semi-free resolution $G$, which has a semi-basis concentrated in degree $u$. Then
\begin{align*}
\tau^{\ge u}M:\qquad\qquad  \cdots\to 0\to H^u\oplus C^u\stackrel{\partial_M^u}{\to} M^{u+1}\stackrel{\partial_M^{u+1}}{\to}
M^{u+2}\stackrel{\partial_M^{u+2}}{\to}\cdots
\end{align*}
is quasi-isomorphic to $M^{\ge t}$. So $\tau^{\ge u}M$ is also a Koszul DG $A$-module and then $\tau^{\ge u}M$ is compact by Proposition \ref{koszulcomp}.

If $H^t(M^{\ge t})=\mathrm{ker}(\partial_{M}^t)\neq 0$, then $M^{\ge t}$ has a minimal semi-free resolution $F_{M^{\ge t}}$, which admits a semi-basis concentrated in degree $t$.
 We have the following short exact sequence
 \begin{align}\label{shortex} 0\to \tau^{\ge t}M \stackrel{\iota}{\to} M^{\ge t}\stackrel{\pi}{\to}  \Sigma^{-t}B^t \to 0 .
 \end{align}
 Let $F_{\k}$ be the minimal semi-free resolution of $\k_A$. Since $A$ is Koszul, $F_{\k}$ admits a semi-basis concentrated in degree $0$.  Acting on
(\ref{shortex}) by $F_{\k}\otimes_{A}-$ gives a new short exact
sequence
$$
0 \to F_{\k}\otimes_{A}\tau^{\ge t}M \stackrel{F_{\k}\otimes_{A}\iota}{\to} F_{\k}\otimes_{A} M^{\ge t} \stackrel{F_{\k}\otimes_{A}\pi}{\to} F_{\k}\otimes_A \Sigma^{-t}B^t
\to 0.
$$
It induces a long exact sequence of cohomologies
\begin{align*}\label{longex}
\cdots \stackrel{\delta^{i-1}}{\to} H^{i}(F_{\k}\otimes_{A}\tau^{\ge t}M)\stackrel{H^{i}(F_{\k}\otimes_{A}\iota)}{\to}
H^{i}(F_{\k}\otimes_{A} M^{\ge t})\stackrel{H^i(F_{\k}\otimes_{A}\pi)}{\to}\\
 H^i(F_{\k}\otimes_A \Sigma^{-t}B^t)\stackrel{\delta^i}{\to}  H^{i+1}(F_{\k}\otimes_{A}\tau^{\ge t}M)\stackrel{H^{i+1}(F_{\k}\otimes_{A}\iota)}{\to} \cdots.
\end{align*}
By the minimality of $F_{\k}$ and $F_{M^{\ge t}}$, we have $$
H(F_{\k}\otimes_{A} M^{\ge t})=H(\k\,{}^L{\otimes}_AM^{\ge t})=H(\k{\otimes}_AF_{M^{\ge t}})=\k\otimes_AF_{M^{\ge t}}
$$
and
$$
H(F_{\k}\otimes_A \Sigma^{-t}B^t)=\Sigma^{-t}(F_{\k}\otimes_AB^t).
$$
Hence $H^i(F_{\k}\otimes_{A} M^{\ge t})=0$ and $H^i(F_{\k}\otimes_A \Sigma^{-t}B^t)=0$ when $i\neq t$. By the definition of connecting homomorphism, it is easy for one to check that $\delta^t=0$. Therefore, $H^i(F_{\k}\otimes_{A}\tau^{\ge t}M)=0$, for any $i\neq t$. If $H^t(F_{\k}\otimes_{A}\tau^{\ge t}M)=0$, then $F_{\k}\otimes_{A}\tau^{\ge t}M$ is quasi-trivial and hence $\tau^{\ge t}M$ is a quasi-trivial DG $A$-module by Lemma \ref{quasitrivial}.  So $\tau^{\ge t}M$ is a zero object in $\mathrm{D}^C(A)$. If $H^t(F_{\k}\otimes_{A}\tau^{\ge t}M)\neq 0$, then
\begin{align*}
H(\k\, {}^L{\otimes}_A\tau^{\ge t}M)=H(F_{\k}\otimes_{A}\tau^{\ge t}M)
\end{align*}
is concentrated in degree $t$. Hence $\tau^{\ge t}M$ is a Koszul DG $A$-module. Then $\tau^{\ge t}M$ is a compact DG $A$-module by Proposition \ref{koszulcomp}.

We have proved the claim that there is some $i\ge b$ such that $\tau^{\ge i}M$ is a quasi-trivial DG $A$-module or a compact DG $A$-module.
On the other hand, the quotient DG $A$-module $M/\tau^{\ge i}M$ is
\begin{align*}
  \qquad\qquad  \cdots \stackrel{\partial_M^{i-3}}{\to} M^{i-2}\stackrel{\partial_M^{i-2}}{\to}
M^{i-1}\stackrel{\partial_M^{i-1}}{\to} B^i\to 0.
\end{align*}
So we have $$H^j(M/\tau^{\ge i}M)=\begin{cases}
0, \quad\quad\quad\quad \text{if}\quad  j< b, \\
H^j(M), \quad\,\,  \text{if}\quad  b\le j\le i-1, \\
0,  \quad\quad\quad\quad \text{if}\quad j\ge i.
\end{cases} $$
Since $M$ is an object in $\mathrm{D}^f(A)$ and $H(A)$ is a Noetherian connected graded algebra, we have $\dim_{\k}H^j(M)<\infty$.
So $M/\tau^{\ge i}M$ is an object in $\mathrm{D^b_{lf}}(A)$.
By \cite[Proposition 7.3]{MW2}, $M/\tau^{\ge i}M$ is a compact DG $A$-module.
The short exact sequence $$0\to \tau^{\ge i}M \to M\to M/\tau^{\ge i}M\to 0$$  implies that $M$ is a compact DG $A$-module.
\end{proof}

 Both
Example \ref{ex1} in the last section and \cite[Example 3.12]{MW2}
 show that the condition that $H(A)$ is a left Noetherian graded
algebra with finite global dimension is much  stronger than the
condition that $A$ is  homologically smooth. A natural question is
whether $\mathrm{D_{fg}}(A)=\mathrm{D}^c(A)$ is still right under the latter
weaker condition. We have the following theorem.

\begin{thm}\label{criteria}
Let $A$ be a Koszul connected cochain DG algebra such that $H(A)$ is a
Noetherian graded algebra with a balanced dualizing complex. Then
$A$ is homologically smooth if and only if
$\mathrm{D}^c(A)=\mathrm{D_{fg}}(A)$.
\end{thm}

\begin{proof}
Obviously,  ${}_A\k$ is an object in $\mathrm{D_{fg}}(A)$. Hence $A$ is homologically smooth if $\mathrm{D}^c(A)=\mathrm{D_{fg}}(A)$. Conversely, if $A$ is homologically smooth, then any object in $\mathrm{D_{fg}}(A)$ is compact by Proposition \ref{fgcomp}. On the other hand, any compact DG $A$-module $M$ admits a minimal semi-free resolution $F_M$, which has a finite semi-free basis. So $F_M$ has a semi-free filtration
$$0=F_M(-1)\subset F_M(0)\subset F_M(1)\subset \cdots \subset F_M(t)=F_M,$$
where each $F_M(i)/F_M(i-1)=A\otimes V_i$ is a DG free $A$-module with $\dim_{\k}V_i<\infty$.
We have a sequence of short exact sequences of DG $A$-modules
\begin{equation}\label{ind}
0\to F_M(i-1)\to F_M(i)\to F_M(i)/F_M(i-1)\to 0 , i=1,2, \cdots, t.
\end{equation}
By (\ref{ind}), we can inductively show that each $F_M(i)$ is an object in $\mathrm{D_{fg}}(A)$. Hence $M\in \mathrm{D_{fg}}(A)$ and $\mathrm{D}^c(A)=\mathrm{D_{fg}}(A)$.
\end{proof}

%%% ------------------------------------------------
\section{A counter example}
Suppose that $A$ is a connected cochain DG algebra such that $H(A)$ is a left
Noetherian graded algebra.  By the
existence of Eilenberg-moore resolution, any object in $\mathrm{D_{fg}}(A)$ is compact if $\mathrm{gl.dim} H(A) <\infty$.  A
natural question is whether the converse is right. In this section, we give counter examples to show that it is generally not true.

\begin{ex}\label{ex1}
 Let $A$ be a connected DG algebra such that $A^{\#} = k\langle x,y\rangle/(xy+yx) $ with $|x|=|y|=1$ and its differential $\partial_A$ is defined by $\partial_A(x) = y^2$ and  $\partial_A(y) = 0$.
\end{ex}

By \cite[Example 3.12]{MW2},  $A$
is a Koszul and homologically smooth DG algebra with
 $$H(A)\cong
 k[\lceil x\rceil^2, \lceil y\rceil]/(\lceil y\rceil^2).$$ Since $H(A)$ is a
Noetherian and Gorenstein graded algebra, $H(A)$ has a balanced dualizing complex. By Theorem \ref{criteria},
we have $\mathrm{D_{fg}}(A)=\mathrm{D}^c(A)$ while $\mathrm{gl.dim}H(A)=+\infty$.

\begin{ex}\label{ex2}
Let $A$ be the connected cochain DG algebra such that
$$ A^{\#}= \k \langle x,y\rangle/\left(\begin{array}{ccc}
                            x^2y-(\xi-1) xyx- \xi yx^2            \\
                             xy^2-(\xi-1) yxy-\xi y^2x           \\
                                                 \end {array}\right) $$
is the graded down-up algebra generated by degree $1$ elements $x,y$, and its differential $\partial_A$ is defined by
 $\partial_{A}(x)=y^2$ and $\partial_{A}(y)=0$, where $\xi$ is a fixed primitive cubic root of unity.
\end{ex}

By \cite[Proposition 6.1]{MHLX}, $A$ is a Koszul Calabi-Yau DG algebra. So $A$ is a Koszul homologically smooth DG algebra. And we have $$H(A)=\frac{\k\langle \lceil xy+yx\rceil, \lceil y\rceil\rangle}{\left(\begin{array}{ccc}
                                  \xi \lceil y\rceil \lceil xy+yx\rceil - \lceil xy+yx\rceil \lceil y\rceil \\
                                   \lceil y^2\rceil \\
                                    \end {array}\right)} $$
by \cite[Proposition 5.5]{MHLX}.  So $H(A)$ is Noetherian and has a balanced dualizing complex since it is a Gorenstein graded algebra.
By Theorem \ref{criteria}, we have $\mathrm{D_{fg}}(A)=\mathrm{D}^c(A)$ while $\mathrm{gl.dim}H(A)=+\infty$.

\begin{ex}\label{ex3}
Let $A$ be the connected cochain DG algebra such that
$$ A^{\#}= \k \langle x,y\rangle/\left(\begin{array}{ccc}
                            x^2y- yx^2            \\
                             xy^2- y^2x           \\
                                                 \end {array}\right) $$
is the graded down-up algebra generated by degree $1$ elements $x,y$, and its differential $\partial_A$ is defined by
 $\partial_{A}(x)=y^2$ and $\partial_{A}(y)=0$.
 \end{ex}
By \cite[Proposition 6.9]{MHLX}, $A$ is a non-Koszul Calabi-Yau DG algebra. So $A$ is a homologically smooth DG algebra.
And we have $$H(A)=\k[\lceil x^2\rceil,\lceil y\rceil,\lceil xy+yx\rceil]/(\lceil y\rceil^2)$$ by \cite[Proposition 5.7]{MHLX}.
So $H(A)$ is Noetherian and has a balanced dualizing complex since it is a Gorenstein graded algebra. By Theorem \ref{koszulcomp}, any Koszul DG module in $\mathrm{D_{fg}}(A)$ is compact. However,  we are unable to judge whether $\mathrm{D_{fg}}(A)=\mathrm{D}^c(A)$ since $A$ is not a Koszul DG algebra.

%%% ------------------------------------------------
\section*{Acknowledgments}
 The first author is
supported by NSFC  (Grant No.11001056),
the
China Postdoctoral
Science Foundation  (Grant Nos.
20090450066 and 201003244),  the
Key Disciplines of Shanghai Municipality (Grant No.S30104) and the Innovation Program of
Shanghai Municipal Education Commission (Grant No.12YZ031).

%%% ------------------------------------------------

\def\refname{References}

\end{document}